\documentclass{ifacconf}
\usepackage{graphicx}      
\usepackage{natbib}        
\usepackage[utf8]{inputenc}

\usepackage{amsmath}
\usepackage{amssymb}

\usepackage{xspace}

\begin{document}
\begin{frontmatter}

\title{Discrete-time Flatness and Linearization along Trajectories} 

\thanks[footnoteinfo]{ This work has been supported by the Austrian Science Fund (FWF) under grant number P 32151.}

\author[First]{Bernd Kolar}
\author[Second]{Johannes Diwold}
\author[Second]{Conrad Gstöttner}
\author[Second]{Markus Schöberl}

\address[First]{Magna Powertrain Engineering Center Steyr GmbH \& Co KG, Steyrer Str. 32, 4300 St. Valentin, Austria \\\ (e-mail: bernd\_kolar@ifac-mail.org)}
\address[Second]{Institute of Automatic Control and Control Systems Technology,	Johannes Kepler University Linz, Altenbergerstraße 66, 4040 Linz, Austria (e-mail: johannes.diwold@jku.at, conrad.gstoettner@jku.at, markus.schoeberl@jku.at)}

\begin{abstract}
The paper studies the relation between a nonlinear time-varying flat discrete-time system and the corresponding linear time-varying systems which are obtained by a linearization along trajectories. It is motivated by the continuous-time case, where it is well-known that the linearization of flat systems along trajectories results in linear time-varying systems which are again flat. Since flatness implies controllability, this constitutes an important verifiable necessary condition for flatness.
In the present contribution, it is shown that this is also true in the discrete-time case: We prove that the linearized system is again flat, and that a possible flat output is given by the linearization of a flat output of the nonlinear system. Analogously, the map that describes the parameterization of the system variables of the linear system by this flat output coincides with the linearization of the corresponding map of the nonlinear system. The results are illustrated by two examples.
\end{abstract}


\end{frontmatter}

\section{Introduction}

The concept of flatness has been introduced in the 1990s by Fliess,
Levine, Martin and Rouchon for nonlinear continuous-time systems,
see e.g. \cite{FliessLevineMartinRouchon:1992}, \cite{FliessLevineMartinRouchon:1995},
or \cite{FliessLevineMartinRouchon:1999}. Since flatness allows an
elegant solution for motion planning problems and a systematic design
of tracking controllers, it is of high practical relevance and belongs
to the most popular nonlinear control concepts. Nevertheless, checking
the flatness of a nonlinear multi-input system is known as a highly
nontrivial problem, for which still no complete systematic solution
in the form of verifiable necessary and sufficient conditions exists
(see e.g. \cite{NicolauRespondek:2016}, \cite{NicolauRespondek:2017},
or \cite{GstottnerKolarSchoberl:2021b} for recent contributions in
this field). For this reason, also necessary conditions for flatness
are of interest to be able to prove at least that a given system is
not flat. One such necessary condition is based on the fact that the
linearization of a flat continuous-time system along a trajectory
yields a linear time-varying system which is again flat and hence
controllable (see e.g. \cite{Rudolph:2021}). Since the latter property
can be checked easily for linear systems, this connection between
a nonlinear system and its linearization constitutes an important
necessary condition for the flatness of continuous-time systems.

The purpose of the present contribution is to investigate the relation
between a flat system and its linearization along a trajectory in
the discrete-time case. Since the linearization of a nonlinear system
along a trajectory leads in general to a linear time-varying system
but the literature has addressed so far only the time-invariant case,
we first need to discuss the concept of discrete-time flatness for
time-varying systems. As proposed in \cite{DiwoldKolarSchoberl:2020},
we consider discrete-time flatness as the existence of a one-to-one
correspondence of the system trajectories to the trajectories of a
trivial system. This leads naturally to a formulation which takes
into account both forward- and backward-shifts of the system variables
as it is also proposed in \cite{GuillotMillerioux:2020}. The point
of view adopted e.g. in \cite{Sira-RamirezAgrawal:2004}, \cite{KaldmaeKotta:2013},
or \cite{KolarKaldmaeSchoberlKottaSchlacher:2016}, where discrete-time
flatness is defined by replacing the time derivatives of the continuous-time
case by forward-shifts, is included as a special case and denoted
within the present paper as forward-flatness.

As our main result, we prove that the linearization of a flat discrete-time
system along a trajectory is again flat, and that a possible flat
output is given by the linearization of a flat output of the nonlinear
system. Furthermore, we show that the corresponding parameterization
of the system variables by the flat output and its shifts coincides
with the linearization of the parameterization of the nonlinear system.
Like in the continuous-time case, this connection between nonlinear
system and linearized system establishes an important necessary condition
for flatness. Even though for discrete-time systems the property of
forward-flatness can be checked efficiently by a generalization of
the test for static feedback linearizability (see \cite{KolarDiwoldSchoberl:2019})
which is based on a certain decomposition property derived in \cite{KolarSchoberlDiwold:2019},
for the more general case including both forward- and backward-shifts
of the system variables a computationally feasible test does not yet
exist.\footnote{An interesting approach can be found in \cite{Kaldmae:2021} but requires
the solution of partial differential equations.} Hence, as we shall illustrate by our second example, the derived
necessary condition is a useful possibility to prove that a given
discrete-time system is not flat.

The paper is organized as follows: First, Section 2 deals with the
concept of discrete-time flatness for time-varying systems. The core
of the paper is then contained in Section 3, which studies the relation
between a flat system and the linear time-varying system obtained
by a linearization along a trajectory. The presented results are illustrated
by two examples in Section 4.

\paragraph*{Notation}

Since we apply differential-geometric concepts, we use index notation
and the Einstein summation convention to keep formulas short and readable.
However, to highlight the summation range especially for double sums,
we also frequently indicate the summation explicitly. For coordinates
that represent forward- or backward-shifts of system variables, we
use a notation with subscripts in brackets. For instance, the $\alpha$-th
forward- or backward-shift of a component $y^{j}$, $j\in\{1,\ldots,m\}$
of a flat output $y$ with $\alpha\in\mathbb{Z}$ is denoted by $y_{[\alpha]}^{j}$,
and $y_{[\alpha]}=(y_{[\alpha]}^{1},\ldots,y_{[\alpha]}^{m})$. Furthermore,
to facilitate the handling of expressions which depend on different
numbers of shifts of different components of a flat output, we use
multi-indices. If $A=(a^{1},\ldots,a^{m})$ is some multi-index, then
$y_{[A]}=(y_{[a^{1}]}^{1},\ldots,y_{[a^{m}]}^{m})$.

\section{Flatness of Time-varying Discrete-time Systems}

In this contribution, we consider nonlinear time-varying discrete-time
systems 
\begin{equation}
x^{i,+}=f^{i}(k,x,u)\,,\quad i=1,\dots,n\label{eq:sys}
\end{equation}
with $\dim(x)=n$, $\dim(u)=m$, and smooth functions $f^{i}(k,x,u)$.
In addition, we assume that the system (\ref{eq:sys}) meets the submersivity
condition
\begin{equation}
\mathrm{rank}(\partial_{(x,u)}f)=n\,,\label{eq:submersivity}
\end{equation}
which is quite common in the discrete-time literature, for all time-steps
$k$.

\begin{figure}

\begin{centering}
\includegraphics[width=8cm]{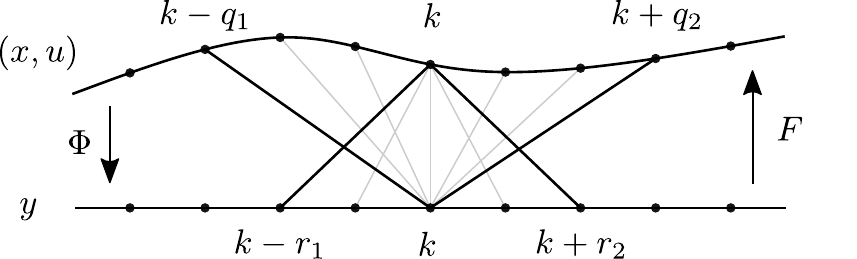}\caption{\label{fig:One-to-one}One-to-one correspondence of the trajectories.}
\par\end{centering}
\end{figure}

As proposed in \cite{DiwoldKolarSchoberl:2020}, where only time-invariant
systems are considered, we call a time-varying discrete-time system
(\ref{eq:sys}) flat if there exists a one-to-one correspondence between
its trajectories $(x(k),u(k))$ and the trajectories $y(k)$ of a
trivial system with $\dim(y)=\dim(u)$. The trajectories of a trivial
system are not restricted by any difference equation and hence completely
free. By one-to-one correspondence, we mean that the values of $x(k)$
and $u(k)$ at a time-step $k$ are determined by an arbitrary but
finite number of future and past values of $y(k)$, i.e., by the trajectory
$y(k)$ in an arbitrarily large but finite time window. Conversely,
the value of $y(k)$ at a time-step $k$ is determined by an arbitrary
but finite number of future and past values of $x(k)$ and $u(k)$.
Consequently, the one-to-one correspondence of the trajectories can
be expressed by maps of the form
\begin{equation}
(x(k),u(k))=F(k,y(k-r_{1}),\ldots,y(k+r_{2}))\label{eq:flat_param_traj}
\end{equation}
and
\begin{equation}
y(k)=\Phi(k,x(k-q_{1}),u(k-q_{1}),\ldots,x(k+q_{2}),u(k+q_{2}))\label{eq:flat_output_traj_xu}
\end{equation}
with suitable integers $r_{1},r_{2},q_{1},q_{2}$ that describe the
length of the corresponding finite time windows, cf. Fig. \ref{fig:One-to-one}.
Since the number of forward- and backward-shifts in (\ref{eq:flat_param_traj})
and (\ref{eq:flat_output_traj_xu}) can of course be different for
the individual components of $y$, $x$, and $u$, we will later use
appropriate multi-indices where it is important.

In the remainder of this section, the framework used in \cite{DiwoldKolarSchoberl:2020}
for the analysis of flat time-invariant discrete-time systems is adapted
to the time-varying case. First, it is important to note that the
representation of a trajectory of the system (\ref{eq:sys}) by both
sequences $x(k)$ and $u(k)$ contains redundancy, as these sequences
are coupled by the system equations (\ref{eq:sys}). By a repeated
application of (\ref{eq:sys}), all forward-shifts $x(k+\alpha)$,
$\alpha\geq1$ of the state variables are determined by $x(k)$ and
the input trajectory $u(k+\alpha)$ for $\alpha\geq0$:
\begin{align*}
x(k+1) & =f(k,x(k),u(k))\\
x(k+2) & =f(k+1,x(k+1),u(k+1))\\
 & \vdots
\end{align*}
In the case $\mathrm{rank}(\partial_{x}f)=n$, the same is also true
for the backward-direction. However, even if the system meets only
the weaker submersivity condition (\ref{eq:submersivity}), there
exist $m$ functions $g(k,x,u)$ such that the $(n+m)\times(n+m)$
Jacobian matrix
\begin{equation}
\begin{bmatrix}\partial_{x}f & \partial_{u}f\\
\partial_{x}g & \partial_{u}g
\end{bmatrix}\label{eq:Jacobian_f_g}
\end{equation}
is regular for all $k$. With such functions, the map
\begin{equation}
\begin{array}{ccc}
x^{+} & = & f(k,x,u)\\
\zeta & = & g(k,x,u)
\end{array}\label{eq:sys_extension}
\end{equation}
is locally invertible for all $k$, and by a repeated application
of its inverse
\begin{equation}
(x,u)=\psi(k,x^{+},\zeta)\label{eq:sys_extension_inverse}
\end{equation}
all backward-shifts $x(k-\beta)$, $u(k-\beta)$, $\beta\geq1$ of
the state- and input variables are determined by $x(k)$ and backward-shifts
$\zeta(k-\beta)$, $\beta\geq1$ of the system variables $\zeta$
defined by (\ref{eq:sys_extension}):
\begin{align*}
(x(k-1),u(k-1)) & =\psi(k-1,x(k),\zeta(k-1))\\
(x(k-2),u(k-2)) & =\psi(k-2,x(k-1),\zeta(k-2))\\
 & \vdots
\end{align*}
Hence, every trajectory of the system (\ref{eq:sys}) is uniquely
determined both in forward- and backward-direction by the values $\ldots,\zeta(k-2),\zeta(k-1),x(k),u(k),u(k+1),\ldots$,
and the map (\ref{eq:flat_output_traj_xu}) can actually be written
as
\begin{equation}
y(k)=\varphi(k,\zeta(k-q_{1}),\ldots,\zeta(k-1),x(k),u(k),\ldots,u(k+q_{2})).\label{eq:flat_output_traj}
\end{equation}
If only a finite time interval is considered, the trajectories of
the system (\ref{eq:sys}) can be identified with points of a finite-dimensional
manifold $\mathcal{\zeta}_{[-l_{\zeta},-1]}\times\mathcal{X}\times\mathcal{U}_{[0,l_{u}]}$
with coordinates $(\zeta_{[-l_{\zeta}]},\ldots,\zeta_{[-1]},x,u,u_{[1]},\dots,u_{[l_{u}]})$
and sufficiently large integers $l_{\zeta}$, $l_{u}$. If
\begin{equation}
h(k,\zeta_{[-l_{1}]},\dots,\zeta_{[-1]},x,u,\dots,u_{[l_{2}]})\label{eq:function_h}
\end{equation}
denotes a function on $\mathbb{Z}\times\mathcal{\zeta}_{[-l_{\zeta},-1]}\times\mathcal{X}\times\mathcal{U}_{[0,l_{u}]}$
which may depend besides the system trajectory also explicitly on
the time-step $k\in\mathbb{Z}$, then its future values can be determined
by a repeated application $\delta^{\beta}$ of the forward-shift operator
$\delta$ defined by the rule
\begin{equation}
\begin{array}{rcll}
k & \rightarrow & k+1\\
\zeta_{[-\beta]}^{j} & \rightarrow & \zeta_{[-\beta+1]}^{j} & \forall\beta\geq2\\
\zeta_{[-1]}^{j} & \rightarrow & g^{j}(k,x,u)\\
x^{i} & \rightarrow & f^{i}(k,x,u)\\
u_{[\alpha]}^{j} & \rightarrow & u_{[\alpha+1]}^{j} & \forall\alpha\geq0\,.
\end{array}\label{eq:forward_shift_operator}
\end{equation}
Likewise, its past values can be determined by a repeated application
$\delta^{-\beta}$ of the backward-shift operator $\delta^{-1}$ defined
by the rule
\begin{equation}
\begin{array}{rcll}
k & \rightarrow & k-1\\
\zeta_{[-\beta]}^{j} & \rightarrow & \zeta_{[-\beta-1]}^{j} & \forall\beta\geq1\\
x^{i} & \rightarrow & \psi_{x}^{i}(k-1,x,\zeta_{[-1]})\\
u^{j} & \rightarrow & \psi_{u}^{j}(k-1,x,\zeta_{[-1]})\\
u_{[\alpha]}^{j} & \rightarrow & u_{[\alpha-1]}^{j} & \forall\alpha\geq1\,,
\end{array}\label{eq:backward_shift_operator}
\end{equation}
where $\psi_{x}$ and $\psi_{u}$ are the corresponding components
of (\ref{eq:sys_extension_inverse}).\footnote{Since we use a finite-dimensional framework, it is important to emphasize
that an application of (\ref{eq:forward_shift_operator}) or (\ref{eq:backward_shift_operator})
is only meaningful if the integers $l_{u}$ and $l_{\zeta}$ are chosen
large enough such that the considered function (\ref{eq:function_h})
does not already depend on $u_{[l_{u}]}$ or $\zeta_{[-l_{\zeta}]}$.
This is assumed throughout the contribution.} With this framework, flatness for nonlinear time-varying discrete-time
systems can be defined as follows. Since flatness is a local concept,
in accordance with the discrete-time literature on static and dynamic
feedback linearization, only a suitable neighborhood of an equilibrium
$(x_{0},u_{0})$ (i.e., $x_{0}=f(k,x_{0},u_{0})$ for all $k$) is
considered, cf. e.g. \cite{NijmeijervanderSchaft:1990} or \cite{Aranda-BricaireMoog:2008}.
However, it is important to emphasize that the concept is still meaningful
even in case the conditions do not hold at the equilibrium point itself
due to a singularity.
\begin{defn}
\label{def:Flatness} The system (\ref{eq:sys}) is said to be flat
around an equilibrium $(x_{0},u_{0})$, if the $n+m$ coordinate functions
$x$ and $u$ can be expressed locally by an $m$-tuple of functions
\begin{equation}
y^{j}=\varphi^{j}(k,\zeta_{[-q_{1}]},\dots,\zeta_{[-1]},x,u,\dots,u_{[q_{2}]})\,,\quad j=1,\ldots,m\label{eq:flat_output}
\end{equation}
and their forward-shifts
\[
\begin{array}{ccl}
y_{[1]} & = & \delta(\varphi(k,\zeta_{[-q_{1}]},\dots,\zeta_{[-1]},x,u,\dots,u_{[q_{2}]}))\\
y_{[2]} & = & \delta^{2}(\varphi(k,\zeta_{[-q_{1}]},\dots,\zeta_{[-1]},x,u,\dots,u_{[q_{2}]}))\\
 & \vdots
\end{array}
\]
up to some finite order. The $m$-tuple (\ref{eq:flat_output}) is
called a flat output.
\end{defn}

The definition ensures the existence of both maps (\ref{eq:flat_param_traj})
and (\ref{eq:flat_output_traj}): The map (\ref{eq:flat_output_traj})
is given by (\ref{eq:flat_output}), and the condition that $x$ and
$u$ can be expressed by (\ref{eq:flat_output}) and its shifts necessitates
the existence of a map\footnote{The multi-index $R=(r^{1},\ldots,r^{m})$ contains the number of forward-shifts
of each component of the flat output which is needed to express $x$
and $u$.}
\begin{equation}
\begin{array}{ccll}
x^{i} & = & F_{x}^{i}(k,y,\dots,y_{[R-1]})\,,\quad & i=1,\dots,n\\
u^{j} & = & F_{u}^{j}(k,y,\dots,y_{[R]})\,, & j=1,\dots,m
\end{array}\label{eq:flat_parameterization}
\end{equation}
which corresponds to (\ref{eq:flat_param_traj}). For notational convenience,
we assume like in \cite{DiwoldKolarSchoberl:2020} that the parameterization
(\ref{eq:flat_parameterization}) of $x$ and $u$ depends only on
forward-shifts of the flat output. This is no restriction, since it
can always be achieved by replacing the components of a flat output
by their highest backward-shifts that occur in (\ref{eq:flat_param_traj}).
The fact that $F_{x}$ in (\ref{eq:flat_parameterization}) is independent
of the highest forward-shifts $y_{[R]}$ that are needed to parameterize
the control inputs $u$ follows from an evaluation of the identity
\begin{multline}
F_{x}^{i}(k+1,y_{[1]},\dots,y_{[R]})\\
=f^{i}(k,F_{x}(k,y,\dots,y_{[R-1]}),F_{u}(k,y,\dots,y_{[R]}))\,,\label{eq:parameterization_identity}
\end{multline}
$i=1,\ldots,n$. This identity reflects the fact that (\ref{eq:flat_param_traj})
maps arbitrary trajectories $y(k)$ of the trivial system to trajectories
$(x(k),u(k))$ of the system (\ref{eq:sys}), which, by definition,
must satisfy the difference equation $x(k+1)=f(k,x(k),u(k))$. Furthermore,
it can be shown in the same way as in the time-invariant case in \cite{DiwoldKolarSchoberl:2020}
that the map (\ref{eq:flat_parameterization}) is unique and that
its Jacobian matrix with respect to the variables $y,y_{[1]},\ldots,y_{[R]}$
has rank $n+m$ for all $k$. As a consequence, the Jacobian matrix
of $F_{x}$ with respect to $y,y_{[1]},\ldots,y_{[R-1]}$ has rank
$n$ for all $k$. This property is essential for trajectory planning
tasks: It ensures that for every initial state $x_{i}$ at an arbitrary
time-step $k_{i}$ and every desired final state $x_{f}$ at a time-step
$k_{f}\geq k_{i}+r$ with $r=\max(r^{1},\ldots,r^{m})$ there exists
a trajectory of the flat output such that the set of equations
\[
\begin{array}{ccl}
x_{i} & = & F_{x}(k_{i},y(k_{i}),y(k_{i}+1),\ldots,y(k_{i}+R-1))\\
x_{f} & = & F_{x}(k_{f},y(k_{f}),y(k_{f}+1),\ldots,y(k_{f}+R-1))
\end{array}
\]
is satisfied identically. Hence, for a flat system (\ref{eq:sys})
it is possible to reach every desired state regardless of the initial
state within $r$ time-steps (locally, where the system is flat).
Accordingly, flat systems are locally reachable.

\section{Flatness of the Linearized System}

A linearization of the system (\ref{eq:sys}) along a trajectory $(x(k),u(k))$
yields a linear time-varying system of the form
\begin{equation}
\Delta x^{i,+}=A_{s}^{i}(k)\Delta x^{s}+B_{j}^{i}(k)\Delta u^{j}\,,\quad i=1,\ldots,n\label{eq:sys_lin_traj}
\end{equation}
with
\[
A_{s}^{i}(k)=\left.\partial_{x^{s}}f^{i}\right|_{x=x(k),u=u(k)}
\]
and
\[
B_{j}^{i}(k)=\left.\partial_{u^{j}}f^{i}\right|_{x=x(k),u=u(k)}\,.
\]
For linear time-varying systems (\ref{eq:sys_lin_traj}), the most
general linear flat output has the form\footnote{Note that if a system can be transformed into Brunovsky normal form
then there also exists a flat output which depends only on the state
variables.}
\begin{multline}
\Delta y^{j}=\sum_{l=1}^{m}\sum_{\beta=1}^{q_{1}}a_{l}^{j,\beta}(k)\Delta\zeta_{[-\beta]}^{l}+\sum_{i=1}^{n}b_{i}^{j}(k)\Delta x^{i}+\\
+\sum_{l=1}^{m}\sum_{\alpha=0}^{q_{2}}c_{l}^{j,\alpha}(k)\Delta u_{[\alpha]}^{l}\,,\quad j=1,\ldots,m\,.\label{eq:flat_output_linear}
\end{multline}
The corresponding parameterization
\begin{equation}
\begin{array}{clc}
\Delta x^{i} & =\sum_{l=1}^{m}\sum_{\alpha^{l}=0}^{r^{l}-1}F_{x,l}^{i,\alpha^{l}}(k)\Delta y_{[\alpha^{l}]}^{l}\,,\quad & i=1,\ldots,n\\
\Delta u^{j} & =\sum_{l=1}^{m}\sum_{\alpha^{l}=0}^{r^{l}}F_{u,l}^{j,\alpha^{l}}(k)\Delta y_{[\alpha^{l}]}^{l}\,,\quad & j=1,\ldots,m
\end{array}\label{eq:flat_parameterization_linear}
\end{equation}
of the system variables by the flat output and its shifts is also
linear. The quantities $\Delta\zeta$, which allow like for the original
nonlinear system a minimal parameterization of the past trajectories,
can be chosen directly as
\begin{equation}
\Delta\zeta^{j}=\left.\partial_{x^{i}}g^{j}\right|_{x=x(k),u=u(k)}\Delta x^{i}+\left.\partial_{u^{l}}g^{j}\right|_{x=x(k),u=u(k)}\Delta u^{l}\,,\label{eq:Delta_zeta}
\end{equation}
$j=1,\ldots,m$, with the functions $g(k,x,u)$ of (\ref{eq:sys_extension}).
Because of the regularity of the Jacobian matrix (\ref{eq:Jacobian_f_g}),
the extended system equations
\[
\begin{array}{ccc}
\Delta x^{i,+} & \!\!=\! & \left.\partial_{x^{s}}f^{i}\right|_{x=x(k),u=u(k)}\Delta x^{s}\!\!+\!\left.\partial_{u^{l}}f^{i}\right|_{x=x(k),u=u(k)}\Delta u^{l}\\
\Delta\zeta^{j} & \!\!=\! & \left.\partial_{x^{s}}g^{j}\right|_{x=x(k),u=u(k)}\Delta x^{s}\!\!+\!\left.\partial_{u^{l}}g^{j}\right|_{x=x(k),u=u(k)}\Delta u^{l}
\end{array}
\]
are obviously invertible with respect to $\Delta x$ and $\Delta u$
for all $k$. Thus, forward- and backward-shifts can be defined analogously
to (\ref{eq:forward_shift_operator}) and (\ref{eq:backward_shift_operator}).

The main objective of this section is to prove that if the system
(\ref{eq:sys}) is locally flat in a neighborhood of the considered
trajectory $(x(k),u(k))$,\footnote{As discussed in Section 2, trajectories $(x(k),u(k))$ can be identified
with points on a manifold $\mathcal{\zeta}_{[-l_{\zeta},-1]}\times\mathcal{X}\times\mathcal{U}_{[0,l_{u}]}$.
Thus, it makes indeed sense to talk of a neighborhood of a trajectory.} then the linearized system (\ref{eq:sys_lin_traj}) is also (globally)
flat. More precisely, we show that the linearization of the flat
output (\ref{eq:flat_output}) is a flat output (\ref{eq:flat_output_linear})
for the linearized system (\ref{eq:sys_lin_traj}), and that the corresponding
parameterization (\ref{eq:flat_parameterization_linear}) of the system
variables $\Delta x$ and $\Delta u$ by the flat output (\ref{eq:flat_output_linear})
is given by the linearization of the nonlinear parameterization (\ref{eq:flat_parameterization})
along the corresponding trajectory of the flat output (\ref{eq:flat_output}).
For the following derivations, it is convenient to consider the linearization
process without immediately inserting a particular trajectory of the
nonlinear system. Instead of (\ref{eq:sys_lin_traj}) we then have
\begin{equation}
\Delta x^{i,+}=\partial_{x^{s}}f^{i}(k,x,u)\Delta x^{s}+\partial_{u^{j}}f^{i}(k,x,u)\Delta u^{j}\,,\label{eq:sys_lin}
\end{equation}
with the remaining variables $x$ and $u$ of the nonlinear system
as placeholders for all its possible trajectories. The advantage of
this approach is that the linearization of an arbitrary function (\ref{eq:function_h})
of the system variables of the nonlinear system (\ref{eq:sys}) can
then formally be written as a Lie derivative along the vector field
\begin{equation}
v_{lin}\!=\sum_{j=1}^{m}\sum_{\beta\geq1}\Delta\zeta_{[-\beta]}^{j}\partial_{\zeta_{[-\beta]}^{j}}\!\!\!+\sum_{i=1}^{n}\Delta x^{i}\partial_{x^{i}}+\sum_{j=1}^{m}\sum_{\alpha\geq0}\Delta u_{[\alpha]}^{j}\partial_{u_{[\alpha]}^{j}},\label{eq:v_lin}
\end{equation}
which considerably simplifies the following calculations. The starting
point is the identity
\begin{equation}
\begin{array}{cll}
x^{i} & =F_{x}^{i}(k,\varphi,\ldots,\delta^{R-1}(\varphi))\,,\quad & i=1,\ldots,n\\
u^{j} & =F_{u}^{j}(k,\varphi,\ldots,\delta^{R}(\varphi))\,, & j=1,\ldots,m
\end{array}\label{eq:identity_nonlinear}
\end{equation}
for the nonlinear system (\ref{eq:sys}), which simply states that
the variables $x$ and $u$ can be expressed by substituting the corresponding
shifts of the flat output (\ref{eq:flat_output}) into the parameterization
(\ref{eq:flat_parameterization}). Computing for both sides of the
identity (\ref{eq:identity_nonlinear}) the Lie derivative along the
vector field (\ref{eq:v_lin}) yields by an application of the chain
rule the identity
\begin{equation}
\begin{array}{clc}
\Delta x^{i} & \!\!=\sum_{l=1}^{m}\sum_{\alpha^{l}=0}^{r^{l}-1}\left(\partial_{y_{[\alpha^{l}]}^{l}}F_{x}^{i}\circ\varphi_{[0,R-1]}\right)\Delta y_{[\alpha^{l}]}^{l}\,, & i=1,\ldots,n\\
\Delta u^{j} & \!\!=\sum_{l=1}^{m}\sum_{\alpha^{l}=0}^{r^{l}}\left(\partial_{y_{[\alpha^{l}]}^{l}}F_{u}^{j}\circ\varphi_{[0,R]}\right)\Delta y_{[\alpha^{l}]}^{l}\,, & j=1,\ldots,m
\end{array}\label{eq:param_lin_initial}
\end{equation}
with
\begin{equation}
\Delta y_{[\alpha]}^{j}=L_{v_{lin}}\delta^{\alpha}(\varphi^{j})\,,\quad\alpha\geq0,\,j=1,\ldots,m\,.\label{eq:delta_y_initial}
\end{equation}
The composition with $\varphi_{[0,R-1]}$ and $\varphi_{[0,R]}$ simply
denotes the substitution of $y,\ldots,y_{[R]}$ by $\varphi,\ldots,\delta^{R}(\varphi)$
as in (\ref{eq:identity_nonlinear}) (note the slight abuse of notation,
since the argument $k$ is not substituted). After substituting a
trajectory of the nonlinear system (\ref{eq:sys}) into the expressions
in brackets, the identity (\ref{eq:param_lin_initial}) would already
look like a parameterization (\ref{eq:flat_parameterization_linear})
of the state- and input variables of the linearized system (\ref{eq:sys_lin_traj})
by a flat output (\ref{eq:flat_output_linear}) and its forward-shifts.
However, it is important to note that $\Delta y_{[\alpha]}^{j}$ denotes
here according to (\ref{eq:delta_y_initial}) only the linearization
of the $\alpha$-th forward-shift of the flat output (\ref{eq:flat_output})
of the nonlinear system. Thus, we must prove that (\ref{eq:delta_y_initial})
coincides with the $\alpha$-th forward-shift of the linearized flat
output
\begin{equation}
\Delta y^{j}=L_{v_{lin}}(\varphi^{j})\,,\quad j=1,\ldots,m\label{eq:linearized_flat_output}
\end{equation}
along trajectories of the linearized system. Since we work with the
linearized equations (\ref{eq:sys_lin}) without a restriction to
a particular trajectory of the nonlinear system, the corresponding
forward-shift operator, which we denote in the following as $\delta_{lin}$,
must shift correctly both the variables $\dots,\Delta\zeta_{[-1]},\Delta x,\Delta u,\Delta u_{[1]},\dots$
of the linearized system as well as the remaining variables of the
nonlinear system. Since the latter serve as placeholders for trajectories
of the nonlinear system, they have to be shifted according to the
rule (\ref{eq:forward_shift_operator}). Thus, the forward-shift operator
$\delta_{lin}$ is defined by the rule
\begin{equation}
\begin{array}{rcll}
k & \rightarrow & k+1\\
\zeta_{[-\beta]}^{j} & \rightarrow & \zeta_{[-\beta+1]}^{j} & \forall\beta\geq2\\
\zeta_{[-1]}^{j} & \rightarrow & g^{j}(k,x,u)\\
x^{i} & \rightarrow & f^{i}(k,x,u)\\
u_{[\alpha]}^{j} & \rightarrow & u_{[\alpha+1]}^{j} & \forall\alpha\geq0\\
\Delta\zeta_{[-\beta]}^{j} & \rightarrow & \Delta\zeta_{[-\beta+1]}^{j} & \forall\beta\geq2\\
\Delta\zeta_{[-1]}^{j} & \rightarrow & \partial_{x^{s}}g^{j}(k,x,u)\Delta x^{s}+\partial_{u^{l}}g^{j}(k,x,u)\Delta u^{l}\!\!\\
\Delta x^{i} & \rightarrow & \partial_{x^{s}}f^{i}(k,x,u)\Delta x^{s}+\partial_{u^{j}}f^{i}(k,x,u)\Delta u^{j}\!\!\\
\Delta u_{[\alpha]}^{j} & \rightarrow & \Delta u_{[\alpha+1]}^{j} & \forall\alpha\geq0\:.
\end{array}\label{eq:delta_lin}
\end{equation}
In the following, we show that for an arbitrary function (\ref{eq:function_h})
we have
\begin{equation}
L_{v_{lin}}(\delta(h))=\delta_{lin}(L_{v_{lin}}(h))\,,\label{eq:commutation_property}
\end{equation}
i.e., shifting along trajectories of the nonlinear system and a subsequent
linearization yields the same result as a linearization and a subsequent
shift along trajectories of the linearized system. If this property
holds for one-fold shifts, then a repeated application immediately
yields the desired result
\begin{equation}
L_{v_{lin}}(\delta^{\alpha}(\varphi))=\delta_{lin}^{\alpha}(L_{v_{lin}}(\varphi))\,,\quad\alpha\geq1\,.\label{eq:commutation_property_flat}
\end{equation}
To prove (\ref{eq:commutation_property}), we simply evaluate both
sides and show that they are equal. Let us start with the right-hand
side: A linearization of the function (\ref{eq:function_h}) yields
\begin{align*}
L_{v_{lin}}(h)= & \sum_{j=1}^{m}\sum_{\beta\geq1}\left(\partial_{\zeta_{[-\beta]}^{j}}h\right)\Delta\zeta_{[-\beta]}^{j}+\sum_{i=1}^{n}\left(\partial_{x^{i}}h\right)\Delta x^{i}+\\
 & +\sum_{j=1}^{m}\sum_{\alpha\geq0}\left(\partial_{u_{[\alpha]}^{j}}h\right)\Delta u_{[\alpha]}^{j}\,,
\end{align*}
and a subsequent shift operation according to (\ref{eq:delta_lin})
results in
\begin{align}
\delta_{lin}(L_{v_{lin}}(h))= & \sum_{j=1}^{m}\sum_{\beta\geq2}\delta(\partial_{\zeta_{[-\beta]}^{j}}h)\Delta\zeta_{[-\beta+1]}^{j}+\nonumber \\
 & +\sum_{j=1}^{m}\delta(\partial_{\zeta_{[-1]}^{j}}h)\delta_{lin}(\Delta\zeta_{[-1]}^{j})+\nonumber \\
 & +\sum_{i=1}^{n}\delta(\partial_{x^{i}}h)\delta_{lin}(\Delta x^{i})+\label{eq:shift_after_lin}\\
 & +\sum_{j=1}^{m}\sum_{\alpha\geq0}\delta(\partial_{u_{[\alpha]}^{j}}h)\Delta u_{[\alpha+1]}^{j}\,.\nonumber 
\end{align}
Note that in order to keep the expression short and facilitate a
comparison, we did not evaluate $\delta_{lin}(\Delta\zeta_{[-1]}^{j})$
and $\delta_{lin}(\Delta x^{i})$. Now let us evaluate the left-hand
side of (\ref{eq:commutation_property}). The forward-shift of (\ref{eq:function_h})
is given by
\begin{multline*}
\delta(h)=\\
h(k+1,\zeta_{[-l_{1}+1]},\dots,\underbrace{g(k,x,u)}_{\zeta},\underbrace{f(k,x,u)}_{x^{+}},u_{[1]},\dots,u_{[l_{2}+1]}),
\end{multline*}
and a subsequent linearization yields
\begin{align}
L_{v_{lin}}(\delta(h))= & \sum_{j=1}^{m}\sum_{\beta\geq1}\left(\partial_{\zeta_{[-\beta]}^{j}}\delta(h)\right)\Delta\zeta_{[-\beta]}^{j}+\nonumber \\
 & +\sum_{j=1}^{m}\left(\left.\partial_{\zeta^{j}}\delta(h)\right|_{\zeta=g}\right)L_{v_{lin}}(g^{j})+\label{eq:lin_after_shift}\\
 & +\sum_{i=1}^{n}\left(\left.\partial_{x^{i,+}}\delta(h)\right|_{x^{+}=f}\right)L_{v_{lin}}(f^{i})+\nonumber \\
 & +\sum_{j=1}^{m}\sum_{\alpha\geq1}\left(\partial_{u_{[\alpha]}^{j}}\delta(h)\right)\Delta u_{[\alpha]}^{j}\,.\nonumber 
\end{align}
By the definition of the forward-shift operator $\delta$ according
to (\ref{eq:forward_shift_operator}), it is straightforward to verify
that
\begin{align*}
\partial_{\zeta_{[-\beta]}^{j}}\delta(h) & =\delta(\partial_{\zeta_{[-\beta-1]}^{j}}h) & \forall\beta\geq1\,\,\\
\left.\partial_{\zeta^{j}}\delta(h)\right|_{\zeta=g} & =\delta(\partial_{\zeta_{[-1]}^{j}}h)\\
\left.\partial_{x^{i,+}}\delta(h)\right|_{x^{+}=f} & =\delta(\partial_{x^{i}}h)\\
\partial_{u_{[\alpha]}^{j}}\delta(h) & =\delta(\partial_{u_{[\alpha-1]}^{j}}h) & \forall\alpha\geq1\,.
\end{align*}
Together with $L_{v_{lin}}(g^{j})=\delta_{lin}(\Delta\zeta_{[-1]}^{j})$
and $L_{v_{lin}}(f^{i})=\delta_{lin}(\Delta x^{i})$ (cf. (\ref{eq:Delta_zeta})
and (\ref{eq:sys_lin})), it can thus be observed that (\ref{eq:shift_after_lin})
and (\ref{eq:lin_after_shift}) are equal, which proves (\ref{eq:commutation_property}).
With (\ref{eq:commutation_property_flat}) it follows then immediately
that the quantities $\Delta y_{[\alpha]}^{j}$ in (\ref{eq:param_lin_initial})
can also be interpreted as forward-shifts of the linearized flat output
(\ref{eq:linearized_flat_output}) along trajectories of the linearized
system. Consequently, after substituting the considered trajectory
of the nonlinear system (\ref{eq:sys}) into the expressions in brackets
of (\ref{eq:param_lin_initial}) as well as into (\ref{eq:linearized_flat_output}),
we have a map (\ref{eq:flat_parameterization_linear}) which allows
to express the state- and input variables $\Delta x$ and $\Delta u$
of the linearized system (\ref{eq:sys_lin_traj}) by a flat output
(\ref{eq:flat_output_linear}) and its forward-shifts. Hence, according
to Definition \ref{def:Flatness}, the linearized system (\ref{eq:sys_lin_traj})
is flat.

As already discussed in Section 2, a particularly important property
of flat systems is the fact that after substituting the parameterization
(\ref{eq:flat_parameterization}) into the system equations (\ref{eq:sys})
the latter are satisfied identically, cf. (\ref{eq:parameterization_identity}).
This can be written formally as
\begin{multline}
\delta_{y}(F_{x}^{i}(k,y,\ldots,y_{[R-1]}))\\
=f^{i}(k,F_{x}(k,y,\ldots,y_{[R-1]}),F_{u}(k,y,\ldots,y_{[R]}))\,,\label{eq:identity_param_nonlinear}
\end{multline}
$i=1,\ldots,n$, with $\delta_{y}$ denoting the forward-shift operator
in $y$-coordinates, which is defined by the rule
\begin{align*}
k & \rightarrow k+1\\
y_{[\alpha]}^{j} & \rightarrow y_{[\alpha+1]}^{j}\,,\quad\forall\alpha\in\mathbb{Z}\,.
\end{align*}
In the following, we show that this property holds indeed also for
the linearized system. More precisely, we show that the linearized
parameterization (\ref{eq:param_lin_initial}) satisfies the linearized
system equations (\ref{eq:sys_lin}) identically (again, it is convenient
to perform the calculations without inserting a particular trajectory
of the nonlinear system (\ref{eq:sys})). First, computing the derivative
of both sides of (\ref{eq:identity_param_nonlinear}) with respect
to $y_{[\alpha^{j}]}^{j}$ for some $j\in\{1,\ldots,m\}$ and $\alpha^{j}\in\{0,\ldots,r^{j}\}$
yields by an application of the chain rule the identity
\[
\partial_{y_{[\alpha^{j}]}^{j}}\delta_{y}(F_{x}^{i})=\left(\partial_{x^{s}}f^{i}\circ F\right)\partial_{y_{[\alpha^{j}]}^{j}}F_{x}^{s}+\left(\partial_{u^{l}}f^{i}\circ F\right)\partial_{y_{[\alpha^{j}]}^{j}}F_{u}^{l}\,.
\]
Since $\delta_{y}$ only substitutes variables, shifting and subsequently
differentiating with respect to $y_{[\alpha^{j}]}^{j}$ is equivalent
to differentiating first with respect to $y_{[\alpha^{j}-1]}^{j}$
and shifting afterwards. Thus, the above identity can be written as
\[
\delta_{y}(\partial_{y_{[\alpha^{j}-1]}^{j}}F_{x}^{i})=\left(\partial_{x^{s}}f^{i}\circ F\right)\partial_{y_{[\alpha^{j}]}^{j}}F_{x}^{s}+\left(\partial_{u^{l}}f^{i}\circ F\right)\partial_{y_{[\alpha^{j}]}^{j}}F_{u}^{l}
\]
for $\alpha^{j}=1,\ldots,r^{j}$ and
\[
0=\left(\partial_{x^{s}}f^{i}\circ F\right)\partial_{y^{j}}F_{x}^{s}+\left(\partial_{u^{l}}f^{i}\circ F\right)\partial_{y^{j}}F_{u}^{l}
\]
for $\alpha^{j}=0$, and a multiplication with $\Delta y_{[\alpha^{j}]}^{j}$
and subsequent summation yields
\begin{multline*}
\sum_{j=1}^{m}\sum_{\alpha^{j}=0}^{r^{j}-1}\delta_{y}\left(\partial_{y_{[\alpha^{j}]}^{j}}F_{x}^{i}\right)\Delta y_{[\alpha^{j}+1]}^{j}\\
=\left(\partial_{x^{s}}f^{i}\circ F\right)\sum_{j=1}^{m}\sum_{\alpha^{j}=0}^{r^{j}-1}\left(\partial_{y_{[\alpha^{j}]}^{j}}F_{x}^{s}\right)\Delta y_{[\alpha^{j}]}^{j}+\\
+\left(\partial_{u^{l}}f^{i}\circ F\right)\sum_{j=1}^{m}\sum_{\alpha^{j}=0}^{r^{j}}\left(\partial_{y_{[\alpha^{j}]}^{j}}F_{u}^{l}\right)\Delta y_{[\alpha^{j}]}^{j}\,.
\end{multline*}
Substituting the flat output $y$ and its shifts by a trajectory
$y(k)$ which corresponds to a considered trajectory of the nonlinear
system (\ref{eq:sys}) finally results in an identity of the form
\begin{multline*}
\sum_{j=1}^{m}\sum_{\alpha^{j}=0}^{r^{j}-1}F_{x,j}^{i,\alpha^{j}}(k+1)\Delta y_{[\alpha^{j}+1]}^{j}\\
=A_{s}^{i}(k)\sum_{j=1}^{m}\sum_{\alpha^{j}=0}^{r^{j}-1}F_{x,j}^{i,\alpha^{j}}(k)\Delta y_{[\alpha^{j}]}^{j}+\\
+B_{l}^{i}(k)\sum_{j=1}^{m}\sum_{\alpha^{j}=0}^{r^{j}}F_{u,j}^{l,\alpha^{j}}(k)\Delta y_{[\alpha^{j}]}^{j}\,,
\end{multline*}
which is the linear equivalent of (\ref{eq:identity_param_nonlinear}).

\section{Examples}

In this section, the derived results are illustrated by two examples.

\subsection{A Flat System}

The purpose of the first example is to show that the linearization
of a flat system along a trajectory is again flat. Since for practically
relevant flat systems like e.g. the gantry crane, the VTOL aircraft,
or the induction motor (see \cite{DiwoldKolarSchoberl:2022} or \cite{DiwoldKolarSchoberl:2020})
the corresponding equations would become rather extensive, for demonstrational
purposes we use the simple academic example
\begin{equation}
\begin{array}{ccl}
x^{1,+} & = & x^{1}+u^{1}\\
x^{2,+} & = & x^{2}+u^{2}\\
x^{3,+} & = & x^{3}+u^{1}u^{2}\,.
\end{array}\label{eq:Product_Example}
\end{equation}
This system corresponds in fact to an exact discretization of the
flat continuous-time system
\begin{equation}
\begin{array}{ccl}
\dot{x}^{1} & = & u^{1}\\
\dot{x}^{2} & = & u^{2}\\
\dot{x}^{3} & = & u^{1}u^{2}
\end{array}\label{eq:Product_Example_cont}
\end{equation}
with a sampling time of $T=1$. With the choice $\zeta^{1}=x^{1}$,
$\zeta^{2}=x^{2}$ for the functions $g(k,x,u)$ of (\ref{eq:sys_extension}),
a flat output of (\ref{eq:Product_Example}) is given by
\begin{equation}
y=(\zeta_{[-1]}^{1},x^{3}-x^{2}(x^{1}-\zeta_{[-1]}^{1}))\,,\label{eq:Product_Example_Flat_Output}
\end{equation}
and the corresponding parameterization of the state- and input variables
(\ref{eq:flat_parameterization}) reads as\footnote{Like the continuous-time system (\ref{eq:Product_Example_cont}),
the system (\ref{eq:Product_Example}) is not flat at equilibrium
points. For an equilibrium the flat output is constant, and the parameterization
(\ref{eq:Product_Example_Param}) becomes singular.}
\begin{align}
x^{1} & =y_{[1]}^{1}\nonumber \\
x^{2} & =\tfrac{y^{2}-y_{[1]}^{2}}{y^{1}-2y_{[1]}^{1}+y_{[2]}^{1}}\nonumber \\
x^{3} & =\tfrac{y^{1}y_{[1]}^{2}-y_{[1]}^{1}(y^{2}+y_{[1]}^{2})+y_{[2]}^{1}y^{2}}{y^{1}-2y_{[1]}^{1}+y_{[2]}^{1}}\label{eq:Product_Example_Param}\\
u^{1} & =y_{[2]}^{1}-y_{[1]}^{1}\nonumber \\
u^{2} & =\tfrac{y^{1}(y_{[1]}^{2}-y_{[2]}^{2})+y_{[1]}^{1}(-y^{2}-y_{[1]}^{2}+2y_{[2]}^{2})}{y^{1}(y_{[1]}^{1}-2y_{[2]}^{1}+y_{[3]}^{1})+y_{[1]}^{1}(-2y_{[1]}^{1}+5y_{[2]}^{1}-2y_{[3]}^{1})+y_{[2]}^{1}(-2y_{[2]}^{1}+y_{[3]}^{1})}\nonumber \\
 & +\tfrac{y_{[2]}^{1}(2y^{2}-y_{[1]}^{2}-y_{[2]}^{2})+y_{[3]}^{1}(-y^{2}+y_{[1]}^{2})}{y^{1}(y_{[1]}^{1}-2y_{[2]}^{1}+y_{[3]}^{1})+y_{[1]}^{1}(-2y_{[1]}^{1}+5y_{[2]}^{1}-2y_{[3]}^{1})+y_{[2]}^{1}(-2y_{[2]}^{1}+y_{[3]}^{1})}\,.\nonumber 
\end{align}
Now let us consider the trajectory
\begin{align}
x^{1}(k) & =\tfrac{1}{2}k(k-1)\nonumber \\
x^{2}(k) & =-\tfrac{1}{2}k(k-1)\nonumber \\
x^{3}(k) & =-\tfrac{1}{6}k(k-1)(2k-1)\label{eq:Product_Example_Trajectory}\\
u^{1}(k) & =k\nonumber \\
u^{2}(k) & =-k\nonumber 
\end{align}
and the corresponding trajectory
\begin{equation}
\begin{array}{cl}
y^{1}(k) & =\tfrac{1}{2}(k-1)(k-2)\\
y^{2}(k) & =\tfrac{1}{6}k(k-1)(k-2)
\end{array}\label{eq:Product_Example_Trajectory_y}
\end{equation}
of the flat output (\ref{eq:Product_Example_Flat_Output}). A linearization
of the system (\ref{eq:Product_Example}) along (\ref{eq:Product_Example_Trajectory})
yields a linear time-varying system (\ref{eq:sys_lin_traj}) with
the matrices
\[
A(k)=\begin{bmatrix}\begin{array}{ccc}
1 & 0 & 0\\
0 & 1 & 0\\
0 & 0 & 1
\end{array}\end{bmatrix}\,,\quad B(k)=\begin{bmatrix}1 & 0\\
0 & 1\\
-k & k
\end{bmatrix}\,,
\]
and a linearization of the flat output (\ref{eq:Product_Example_Flat_Output})
along (\ref{eq:Product_Example_Trajectory}) yields
\begin{equation}
\Delta y=(\Delta\zeta_{[-1]}^{1},\Delta x^{3}-(k-1)\Delta x^{2}+\tfrac{1}{2}k(k-1)(\Delta x^{1}-\Delta\zeta_{[-1]}^{1}))\label{eq:Product_Example_Flat_Output_Linearized}
\end{equation}
with $\Delta\zeta^{1}=\Delta x^{1}$ and $\Delta\zeta^{2}=\Delta x^{2}$
according to (\ref{eq:Delta_zeta}). With a computer algebra program,
it is easy to verify that all state variables $\Delta x^{1},\Delta x^{2},\Delta x^{3}$
and input variables $\Delta u^{1},\Delta u^{2}$ of the linearized
system (\ref{eq:sys_lin_traj}) can be expressed by (\ref{eq:Product_Example_Flat_Output_Linearized})
and its forward-shifts. Furthermore, the corresponding map (\ref{eq:flat_parameterization_linear})
coincides indeed with the linearization
\begin{align*}
\Delta x^{1} & =\Delta y_{[1]}^{1}\\
\Delta x^{2} & =\tfrac{1}{2}k(k-1)(\Delta y^{1}-2\Delta y_{[1]}^{1}+\Delta y_{[2]}^{1})+\Delta y^{2}-\Delta y_{[1]}^{2}\\
\Delta x^{3} & =\tfrac{1}{2}k(k-1)(k\Delta y^{1}-(2k-1)\Delta y_{[1]}^{1}+(k-1)\Delta y_{[2]}^{1})\\
 & +k\Delta y^{2}+(1-k)\Delta y_{[1]}^{2}\\
\Delta u^{1} & =-\Delta y_{[1]}^{1}+\Delta y_{[2]}^{1}\\
\Delta u^{2} & =\tfrac{1}{2}k((1-k)\Delta y^{1}+(3k-1)\Delta y_{[1]}^{1}-(3k+1)\Delta y_{[2]}^{1})\\
 & +\tfrac{1}{2}k(k+1)\Delta y_{[3]}^{1}-\Delta y^{2}+2\Delta y_{[1]}^{2}-\Delta y_{[2]}^{2}
\end{align*}
of the map (\ref{eq:Product_Example_Param}) along the trajectory
(\ref{eq:Product_Example_Trajectory_y}).

\subsection{A Non-Flat System}

As a second example, let us consider the system
\begin{equation}
\begin{array}{ccl}
x^{1,+} & = & -\sin(x^{1}-x^{3})+u^{2}\\
x^{2,+} & = & (1-\sin(x^{1}-x^{3}))u^{1}\\
x^{3,+} & = & u^{2}\,.
\end{array}\label{eq:Non_Flat_Sys}
\end{equation}
A linearization along an arbitrary trajectory $(x(k),u(k))$ results
in a linear time-varying system (\ref{eq:sys_lin_traj}) with\[ A(k)=\!\left.\begin{bmatrix}\begin{array}{ccc} -\cos(x^{1}-x^{3}) & 0 & \cos(x^{1}-x^{3})\\ -\cos(x^{1}-x^{3})u^{1} & 0 & \cos(x^{1}-x^{3})u^{1}\\ 0 & 0 & 0 \end{array}\end{bmatrix}\right|_{\begin{array}{l} {\scriptstyle x=x(k),}\\ {\scriptstyle u=u(k)} \end{array}} \] and \[ B(k)=\!\left.\begin{bmatrix}0 & 1\\ 1-\sin(x^{1}-x^{3}) & 0\\ 0 & 1 \end{bmatrix}\right|_{\begin{array}{l} {\scriptstyle x=x(k),}\\ {\scriptstyle u=u(k)} \end{array}}. \]Since
$A(k)B(k)=0$ for all $k$ independently of the chosen trajectory,
it can be observed immediately that the linearized system (\ref{eq:sys_lin_traj})
is not reachable (conditions for the reachability of linear time-varying
discrete-time systems can be found e.g. in \cite{Weiss:1972}). Thus,
the linearized system cannot be flat, and because of the connection
between the flatness of a nonlinear system and its linearization discussed
in Section 3, the nonlinear system (\ref{eq:Non_Flat_Sys}) cannot
be flat either. This result can also be obtained in an alternative
way by showing e.g. with the method discussed in \cite{Aranda-BricaireKottaMoog:1996}
that the considered nonlinear system (\ref{eq:Non_Flat_Sys}) itself
is also not reachable, and hence clearly not flat.

\section{Conclusion}

We have shown that -- like in the continuous-time case -- the linearization
of a flat discrete-time system yields a linear time-varying system
which is again flat. Since flatness implies reachability (and consequently
also controllability), this property constitutes a useful necessary
condition for flatness. Moreover, we have shown that a possible flat
output can be obtained by a linearization of a flat output of the
nonlinear system, and that the corresponding parameterization of the
system variables of the linearized system coincides with the linearization
of the parameterization of the system variables of the nonlinear system.

\bibliography{Bibliography_Bernd}

\end{document}